\documentclass[11pt]{amsart}
\usepackage{amssymb}

\usepackage{palatino}
\input amssym.def

\usepackage{amsmath, amsfonts}
\usepackage{amssymb}
\usepackage{amscd}
\usepackage[mathscr]{eucal}
\usepackage{palatino}
\newfont{\cyrr}{wncyr10}

\newcommand{\Z}{{\mathbb Z}}
\newcommand{\Q}{{\mathbb Q}}

\newcommand{\N}{{\mathbb N}}
\newcommand{\gH}{{\goth H}}
\newcommand{\gF}{{\goth F}}
\newcommand{\SL}{{\rm SL}}
\newcommand{\bQ}{\overline{\Q}}
\newcommand{\bZ}{\overline{\Z}}

\newtheorem{thm}{Theorem}

\newtheorem{cor}[thm]{Corollary}
\newtheorem{prop}[thm]{Proposition}
\newtheorem{rmk}{Remark}[section] 

\newtheorem{conj}[thm]{Conjecture}

\begin{document}

\title[Zeros of weakly holomorphic modular forms]{On the 
zeros of weakly holomorphic modular forms}

\author{Sanoli Gun and Biswajyoti Saha}
\address{Sanoli Gun and Biswajyoti Saha  $\phantom{mmmmmmmmmm
mmmmmmmmmmmmmmmmmmmmmmmmmmmmmmmmmmm}$\\
Institute of Mathematical Sciences, C.I.T Campus, Taramani, 
Chennai, 600 113, India}
\email{sanoli@imsc.res.in, biswajyoti@imsc.res.in}



\maketitle   

\begin{abstract}
In this article, we study the nature of zeros of weakly 
holomorphic modular forms. In particular, we prove
results about transcendental zeros of
modular forms of higher levels and for certain Fricke groups which
extend a work of Kohnen (see \cite{WK}).  
Furthermore, we investigate the algebraic 
independence of values of weakly holomorphic modular forms. 
\end{abstract}

\section{Introduction}

Throughout the paper, let 
$\gH := \{ z ~ |~  \Im(z) > 0 \}$ be the upper half-plane,
$k$ be an even natural number and 
$\Gamma:= \SL_2(\Z)$ be the full modular
group.  Further, a CM point is an element of 
$\goth H$ lying in an imaginary quadratic field.  
A holomorphic function $f$ on $\gH$ is called a weakly 
holomorphic modular form of weight $k$ for  $\Gamma$ if 
for any $\gamma = \big({ a\ b \atop c\ d }\big) \in \Gamma$ 
$$
f \left(\frac{az+ b}{cz+d} \right) 
= (cz + d)^{k} f(z)  \text{  for all } z \in \gH
$$
and $f$ has an expansion of the form
\begin{equation}\label{rel}
f(z) =  \sum_{n \ge n_0} a(n) e^{2\pi i n z}, 
\text{  where } z \in \gH.
\end{equation}
Moreover, if $n_0 = 0$ (resp. $n_0 = 1$) 
in equation \eqref{rel}, then we call $f$ a 
modular form (resp. cusp form) of weight $k$ for $\Gamma$. 
From now on, we will denote weakly holomorphic modular forms 
by WH modular forms. Also denote $e^{2\pi i z}$ 
by $q$ where $z \in \gH$.
Some well known examples of WH modular forms are
Eisenstein series of weight $k \ge 4$ defined by 
$$
E_k(z) := 1 - \frac{2k}{B_k} 
\sum_{n=1}^{\infty} \sigma_{k-1}(n) q^n, 
$$
where $\sigma_k(n) = \sum_{d|n} d^k$, $B_k$'s are 
Bernoulli numbers,
the normalized weight $12$ Delta function 
$$
\Delta (z) :=  q \prod_{n \ge 1} ( 1 -  q^n)^{24}
$$  
and the weight zero $j$-function given by
$$
j(z) := 1728 ~~\frac{E_4(z)^3}{E_4(z)^3 - E_6(z)^2}~.
$$
For other examples of families of WH modular
forms, see \cite{{AKN}, {DJ}, {KZ}}.

In this article, we study the nature of zeros 
of WH modular forms. From now on, by zeros of 
$f$, we mean inequivalent zeros of $f$, 
i.e. we investigate zeros of $f$ in $\gF$, where
$$
\gF : = \left\{  z \in \gH ~|~  
|z| \ge 1, -1/2 \le  \Re(z) \le 0 \right\} \cup 
\left\{  z \in \gH ~|~  |z| >  1,   0 <  \Re(z) <  1/2 \right\} 
$$
is the standard fundamental domain. 
We will restrict our attention to forms whose 
Fourier coefficients $a(n)$ in 
equation \eqref{rel} are in a subfield $\rm F$ of $\bQ$.
These WH modular forms (resp. modular forms and cusp forms) 
of weight $k$ for $\Gamma$ form a vector space over $\rm F$ denoted
by $M_k^{!}(\rm F)$ (resp. by $M_k(\rm F)$ and $S_k(\rm F)$).

The paper consists of theorems of three different types. 
First we discuss about the existence of at 
least one transcendental 
zero of certain WH modular forms and use this
information to provide a criterion to study
the nature of Fourier coefficients of WH modular
forms.

Next we show
that most of the zeros of certain families
of WH modular forms for level one
and higher level and for certain Fricke
groups (defined later) are transcendental. 
In particular, our results give 
a generalization of a theorem of Kohnen \cite{WK}.
An essential ingredient in the proof
of these theorems is the {\it $q$-expansion principle}
due to Deligne and Rapoport  \cite{DR}
(see Theorem 3.9, p. 304).

In the final set of theorems, we study
the values of WH modular forms at certain
algebraic and transcendental numbers.

The motivation to prove these results
come from the classical results
of Chowla-Selberg \cite{CS}, Chudnovsky \cite{{GV},{GV1}}, 
Nesterenko \cite{YN}, Ramachandra \cite{KR}
and Schneider \cite{TS} and the recent works 
of the first author with Murty
and Rath \cite{GMR}, Kanou \cite{NK} and Kohnen \cite{WK}.

\smallskip

\section{Statement of theorems}

\smallskip

In this section, we list the theorems we prove
in the final section.

\begin{thm}\label{five}
Let $f(z)= \sum_{n \ge n_0} a(n)q^n \in  M_k^{!}(\bQ)$ with
$a(n_0) \in \bZ \setminus \{ 0 \}$ and not all $a(n) \in \bZ$.
Then $f$ has at least one transcendental zero. 
\end{thm}

\begin{rmk}\rm
Let $\beta \in \gH$ be a transcendental number 
such that $j(\beta) \in \bQ \setminus \bZ$.  
Examples of infinitely many such 
transcendental numbers are given in \cite{SG}.
For $z \in \gH$, define a weakly 
holomorphic modular form $f$ of weight $k$ 
such that
$$
f(z) = \Delta^{k/12}(z)\left(j(z) - j(\alpha_1)\right) 
\cdots \left(j(z) - j(\alpha_n)\right)\left(j(z) - j(\beta)\right),
$$ 
where $\alpha_1, \cdots, \alpha_n \in \gH$ are CM points. Then
$f$ satisfies the assumptions of Theorem \ref{five}. Note that
$f$ has exactly one transcendental zero and
hence Theorem \ref{five} is the best possible. 
\end{rmk}

\begin{rmk}\rm
Below we construct examples to show that both the
hypothesis in Theorem \ref{five} are necessary.
Given a finite set of CM points $\alpha_1, \cdots, \alpha_n$ and
$z \in \gH$, we can define a weakly 
holomorphic modular form $f$ of weight $k$
such that
\begin{equation}\label{exam}
f(z) = c~\Delta^{k/12}(z)\left(j(z) - j(\alpha_1)\right) 
\cdots \left(j(z) - j(\alpha_n)\right),
\text{  where } c \in \bQ \setminus \bZ.
\end{equation}
Then $f$ has Fourier coefficients in
$\bQ \setminus \bZ$ and only algebraic zeros. 
By taking $g(z) := f(z)/c$, where $f$ is as in equation
\ref{exam}, we see that $g$ has all coefficients
in $\bZ$ but has only algebraic zeros.
\end{rmk}

Define 
$$
 M_k^{!}(\bZ):=  \left\{ f(z) = \sum_{n \ge n_0} 
 a(n)q^n \in M_k^{!}(\bQ)
~ \mid~   a(n) \in \bZ   ~~~~\forall n \right\}.
$$

Using Theorem \ref{five}, we give a criterion for an 
$f \in M_k^{!}(\bQ)$ 
to have algebraic integer Fourier coefficients.

\begin{cor}\label{six}
Suppose that 
$f(z)= \sum_{n \ge n_0} a(n)q^n$ 
is a WH modular form of weight $k$, 
where $a(n_0)$ is a non-zero algebraic integer. 
Further suppose that all zeros of $f$ are algebraic. 
Then $f \in M_k^{!}(\bQ)$
if and only if  $f \in M_k^{!}(\bZ)$.
\end{cor}

Our next set of theorems show that
information about location of zeros
can be useful to determine the
nature of zeros. These theorems
can be thought of as generalizations
of a theorem of Kohnen (see \cite{WK}).

\begin{thm}\label{seven}
Let $f \in  M_k^{!}(\Q)$ be a weakly 
holomorphic modular form with 
all its zeros on $A$, where 
$$
A:= \{e^{i\theta}\mid \pi/2\leq\theta\leq 2\pi/3\}.
$$
Then all zeros of $f$ other than the possible zeros  
at $i$ and $\rho = e^{2i\pi/3}$ are transcendental.
\end{thm}

\begin{rmk}\rm\label{imp}
Rankin and Swinnerton-Dyer \cite{RS} showed that the
Eisenstein series $E_k$ for even $k \ge 4$ have all their 
zeros on the arc $A$.  See \cite{{AKN}, {DJ}, {KZ}}
for other examples of WH modular 
forms having similar properties. Contrary to these
forms, zeros of Hecke eigen cusp forms are
equidistributed in $\gF$ with respect to the 
hyperbolic measure {\rm(}see \cite{ZR} and \cite{HS}
for further details{\rm)}.
\end{rmk}

Using Theorem \ref{seven}, we now prove a higher 
level analogue. For a natural number $N$, let
$$
\Gamma_0(N) : = 
\left\{  \left( { a\ \ b \atop c\ \ d } 
\right) \in \Gamma 
~\mid~  c \equiv 0\!\!\!\pmod{N}  \right\}.
$$
A holomorphic function $f$ on $\gH$ is 
called a WH modular form of weight $k$ 
and level $N$ if for any 
$\gamma = \big({ a\ b \atop c\ d }\big) \in \Gamma_0(N)$, 
$$
f \left(\frac{az+ b}{cz+d} \right) 
= (cz + d)^{k} f(z)  \text{  for all } z \in \gH
$$
and $f$ is meromorphic at all its cusps. Further, if $f$
is holomorphic at all cusps, we call $f$ a modular form.
We will consider inequivalent zeros of such
WH modular forms in the fundamental domain
$$
\gF^{(N)} := \cup_{ \gamma \in \Gamma / \Gamma_0(N)} \gF.
$$
In this set-up, we prove the following theorem.
\begin{thm}\label{ate}
Let $p$ be a prime and $f$ be a WH modular form 
of weight $k$ and level $p$ having rational Fourier coefficients 
at the cusp $i \infty$. Suppose that all zeros of $f$ lie on 
$$
B:=\cup_{ \gamma \in \Gamma / \Gamma_0(p)} \gamma A,
$$
where $A$ is as in Theorem \ref{seven}.
Let
$$
C := \left\{ \frac{i-n}{n^2+1} ~\mid~  0 \leq n < p \right\} \cup
\left\{\frac{\frac{\sqrt 3}{2}i-n + \frac{1}{2}}{n^2-n+1}  
~\mid~  0\leq n <p \right\}.
$$
Then all zeros of $f$ except the possible 
CM zeros lying in $C$ are transcendental.
\end{thm}

Next we study transcendental zeros of Eisenstein series
for certain Fricke groups $\Gamma_0^*(p)$, where
$$
\Gamma_0^*(p) : = \Gamma_0(p) \cup \Gamma_0(p)W_p, \phantom{m}
W_p:= \begin{pmatrix}
         0       &   -1/\sqrt{p} \\
      \sqrt{p}   &       0  
      \end{pmatrix}
$$
and $p$ is a prime. Here we will study zeros 
in the fundamental domain 
$$
\gF_p := \left\{ |z| \geq 1/{\sqrt{p}}, ~~-1/2 \leq \Re(z) \leq 0 \right\}
\cup \left\{ |z| > 1/{\sqrt{p}}, ~~ 0 < \Re(z) < 1/2\right\}.
$$

One can define modular forms 
on Fricke groups analogously to modular forms for 
congruence subgroups. The Eisenstein series for $\Gamma_0^*(p)$
is defined by
$$
E_{k,p}(z):= \frac{1}{p^{k/2} + 1} \left( p^{k/2} E_k(pz) + E_k(z) \right),
$$ 
where $E_k(z)$ is the Eisenstein series for $\Gamma$.
For $p = 2, 3$, the location of zeros of $E_{k,p}(z)$ were studied by
Miezaki, Nozaki and Shigezumi (see \cite{MNS}). 
They showed that the zeros of $E_{k,p}(z)$ lie on the arc
$$
A_p :=  \left\{ z \in \gH ~|~ -1/2 \leq \Re(z) \leq 0, 
~ |z|= 1/\sqrt{p} \right\}.
$$

The nature of the zeros of the Eisenstein series $E_k$ for the full modular
group $\Gamma$ was first studied by Kanou \cite{NK}. 
He showed that $E_k$ for even $k \ge 16$  
has at least one transcendental zero in $\gF$.  Soon after, 
Kohnen \cite{WK} proved that any zero of $E_k$ in $\gF$  
different from $i$ or $\rho$ is necessarily 
transcendental. See also \cite{{SG}} for 
similar results. Here we prove:

\begin{thm}\label{eisen}
Let $E_{k,2}(z)$ be the Eisenstein series of weight $k$ for the
Fricke group $\Gamma_0^*(2)$. Then all zeros of $E_{k,2}(z)$ 
other than the possible zeros at
$$
\frac{i}{\sqrt{2}}, 
\phantom{m} \frac{-1 + i \sqrt{7}}{4},
\phantom{m} \frac{-1 + i}{2}
$$
are transcendental. 
\end{thm}

Analogously, we have the following theorem for $E_{k,3}(z)$.

\begin{thm}\label{eisen1}
Let $E_{k,3}(z)$ be the Eisenstein series of weight $k$ for the
Fricke group $\Gamma_0^*(3)$. All zeros of $E_{k,3}(z)$ 
are transcendental other than the 
following possible CM zeros
$$
\frac{i}{\sqrt{3}},
\phantom{m} \frac{-1 + i \sqrt{11}}{6}, 
\phantom{m} \frac{-1 + i\sqrt{2}}{3}, 
\phantom{m} \frac{-3 + i \sqrt{3}}{6}.
$$
\end{thm}

The notion of real zeros of modular forms was introduced by 
Ghosh and Sarnak  in \cite{GS} and can be
extended to WH modular forms. 
A zero $z_0$ of a WH modular form $f$ 
is called real if it lies on the arc $A$ or it lies on the 
vertical line passing through $\rho$ or 
the vertical line passing through $i$. In 
this context, we have the following
theorems.

\begin{thm}\label{eight}
Let $f \in  M_k^{!}(\Q)$ be a weakly 
holomorphic modular form with 
all its zeros on $L$, where
$$
L:= -\frac{1}{2} + i t \phantom{m}\text{with } 
t \ge \frac{\sqrt{3}}{2}.
$$
Then any CM zero of $f$ (if exists) is necessarily of the form  
$$
-\frac{1}{2} + i\frac{\sqrt{|D|}}{2a}, \phantom{m}
\text{    where } 1 \le 
a \le \left[\sqrt{\frac{|D|}{3}}\right], ~~a \in \N.
$$
Here $D$ is a discriminant of an imaginary
quadratic field which is necessarily 
congruent to $1\pmod{4}$. 
\end{thm}

\begin{rmk}\rm\label{nine}
Let $f$ be a WH modular form 
as in Theorem \ref{eight}.
Note that the above theorem helps us to 
calculate the possible CM zeros with bounded 
imaginary parts of such a form. In particular 
if we consider the collection
$$
M_L :=  \{  f \in  M_k^{!}(\Q) ~| ~ 
f(\alpha)= 0 \implies \alpha \in L \},
$$
then all zeros of $f \in M_L$ on the line segment
$$ 
L_2:= \{ z \in L ~|~ \Im(z) < 2 \}
$$
are transcendental except for 
$$
\rho, ~~~ -\frac{1}{2} + i \frac{\sqrt{7}}{2}, ~~~
 -\frac{1}{2} + i \frac{\sqrt{11}}{2}, ~~~
 -\frac{1}{2} + i \frac{\sqrt{15}}{2}, ~~~
 -\frac{1}{2} + i \frac{\sqrt{15}}{4}.
$$
\end{rmk}

\pagebreak

\begin{thm}\label{ten}
Let $f \in  M_k^{!}(\Q)$ be a WH modular form with 
all its zeros on 
$$
R :=  it  \phantom{m} \text{with } t \ge 1.
$$
Then any CM zero of $f$ (if exists) is of the form  
$$
i\frac{\sqrt{|D|}}{2a} \phantom{m}
\text{with  } 1 \le a \le 
\left[\frac{\sqrt{|D|}}{2}~\right], ~~a \in \N.
$$
Here $D$ is a discriminant of an imaginary quadratic field
and necessarily $D \equiv 0 \pmod{4}$.
\end{thm}

\begin{rmk}\rm\label{ele}
As before, by applying Theorem \ref{ten} we can calculate possible 
CM zeros with bounded imaginary part for forms of the above type. 
In particular, if we consider the collection
$$
M_R :=  \{ f \in  M_k^{!}(\Q) ~| ~ f(\alpha)= 0 \implies 
\alpha \in R \},
$$
then all zeros of $f \in M_R$ on the line segment
$$ 
R_2:= \{ z \in R ~|~ \Im(z) < 2 \}
$$
other than the zeros at
$$
i, \phantom{m} i\sqrt{2} \phantom{m} \text{and}
\phantom{m} i \sqrt{3}
$$
are transcendental. 
\end{rmk}

\begin{rmk}\rm
Let $g_k$ be any family of WH modular forms having
all their zeros on the arc $A$. See Remark~\ref{imp}
for examples of such families. For $z \in \gH$, let 
$$
f_k^L(z) := g_k(\gamma_L^{-1}z) 
\phantom{m} \text{ and   } \phantom{m}
f_k^R(z) := g_k(\gamma_R^{-1}z), 
$$
where $\gamma_L = \big( { - 2 \ \ -1 \atop  
\  1 \ \ \ -1 } \big)$
and  $\gamma_R = \big( { 1 \ \ -1 \atop  1 \ \ \ \ 1 } \big)$.
Then $f_k^L$ (resp. $f_k^R$) have all their zeros on
$L$ (resp. on $R$).
\end{rmk}

Finally we discuss the values of 
WH modular forms. This can be done
by suitably modifying the arguments followed in 
the work of the first author, 
Murty and Rath \cite{GMR}. 
For the sake of completion, we indicate 
the relevant modifications
and sketch the proofs. 
In order to do so, let us define an 
equivalence relation on $M^!(\rm F)$, the graded 
ring (graded by the weight $k$) 
of WH modular forms with Fourier coefficients in 
a subfield $\rm F$ of $\bQ$
as follows.  Two such WH modular forms $f$ and 
$g$ are called equivalent, 
denoted by  $f \sim g$, if there exists 
natural numbers $k_1, k_2$ such that 
$f^{k_2} = \lambda g^{k_1}$
with $\lambda \in {\rm F}^*$. Otherwise 
we call $f$ not equivalent to $g$
and write $f \not\sim g$.  The restrictions 
of this relation to the
graded rings of modular forms and 
cusp forms with Fourier coefficients in $\rm F$ 
give rise to equivalence relations on those rings.

\begin{thm}\label{one}
Let $f$ be any non-zero element of $M_k^{!}(\bQ)$. Suppose 
that $f \not\sim \Delta$ and $\alpha \in {\goth H}$ is algebraic.
Then $f^{12}(\alpha)/\Delta^k(\alpha)$ 
is algebraic if and only if $\alpha$
is a {\rm CM} point.
\end{thm}

\begin{rmk}\rm
It was noticed in \cite{GMR} that if $\alpha$ is a CM point, 
then $\Delta(\alpha)$ is transcendental 
by a theorem of Schneider (see \cite{TS}).
Further it was shown that if $\alpha$ is transcendental 
with $j(\alpha) \in \bQ$, then also $\Delta(\alpha)$ is transcendental. 
Moreover, if $\alpha \in \goth H$ is non-CM algebraic, a conjecture 
of Nesterenko (see \cite{NP}, page 31) will imply the 
transcendence of $\Delta(\alpha)$.  Using the Open Mapping Theorem,
we see that $\Delta(\alpha)$ can take algebraic values. It is then 
clear from the above discussion that in this case 
$\alpha$ is conjecturally transcendental 
and $j(\alpha)$ is transcendental. 
\end{rmk}
 
Let $\alpha \in \gH$.  Now if  $\alpha \in \bQ$ or $j(\alpha) \in \bQ$,
we can deduce the nature of $f(\alpha)$ except in one case. 
More precisely, we have the following theorems.

\begin{thm}\label{three}
Let $f \in M_k^{!}(\bQ)$ be non-zero and $\alpha \in \gH$ be 
such that $j(\alpha) \in \bQ$. 
Then either $f(\alpha) = 0$ or $f(\alpha)$ is algebraically 
independent with $e^{2\pi i \alpha}$.
\end{thm}

\begin{thm}\label{four}
Let $\alpha \in {\goth H}$ be a non-CM algebraic number.  
Also, let 
$$
{\rm S}_\alpha : = 
\left\{ f \in M^{!}(\bQ) ~|~ f \ne 0 \text{  and } 
f(\alpha) \in \bQ  \right\} / \sim.
$$ 
Then ${\rm S}_\alpha$ has at most one element.
\end{thm}

\begin{rmk}\rm
We note that a conjecture of Nesterenko (discussed in Section 3) will 
imply that $\rm S_{\alpha}$ is empty as was noticed earlier in \cite{GMR}
in the case of modular forms.
\end{rmk}

Finally, we have the following general theorem about the nature of zeros 
of WH modular forms. 
 
\begin{thm}\label{two}
Let $f \in M_k^{!}(\bQ)$ be non-zero. Any zero of $f$ is 
either {\rm CM} or transcendental.
\end{thm}

\smallskip

\section{Preliminaries}

\smallskip

We begin by fixing some notations and recalling some results. 
As before, an element $\alpha \in \gH$ lying in an imaginary quadratic 
field $K = \Q(\sqrt{D})$ will be called a CM point.  
It is known, from classical theory of complex multiplication 
that if $\alpha \in \gH$ is a CM point, then $j(\alpha)$ is an algebraic integer,
lying in the Hilbert class field of $\overline{\Q(\alpha)}$, the Galois closure of
$\Q(\alpha)$.  
See Chapter 3 of \cite{DC} and Chapter 10 of \cite{SL}
for further details on the theory of complex multiplication.
Moreover, from Theorem 5 of \cite{SL}, one can 
get the following proposition.

\begin{prop}\label{cm}
Let $\alpha_0 \in \gH$ be a CM point and also let $|D|$ be the absolute value
of the discriminant of the lattice $\Lambda := (\alpha_0, 1)$.  Then there is an
element $\sigma \in {\rm Gal}\left(\overline{\Q(\sqrt{D})} / \Q(\sqrt{D})  \right)$
such that $\sigma(j(\alpha_0)) = j(\alpha_1)$, where
$$
\alpha_1:=  \begin{cases}
      \frac{i\sqrt{|D|}}{2} &\text{if $D \equiv 0\!\!\!\pmod 4$,}\\
      \frac{-1+ i\sqrt{|D|}}{2} &\text{if $D \equiv 1\!\!\!\pmod 4$.}
      \end{cases}
$$
\end{prop}

On the other hand, for algebraic points in the upper half plane,
Schneider \cite{TS} proved the following result.

\begin{thm}{\rm (Schneider)}
For $\alpha \in  \gH$, if $\alpha$ and $j(\alpha)$ are algebraic, then 
$\alpha$ is CM.
\end{thm}

The above theorem of Schneider along with the following result
of Nesterenko \cite{YN} play an important role in order to understand
the nature of values of modular forms as shown in \cite{GMR}.

\begin{thm}{\rm (Nesterenko)}
Let $\alpha \in \gH$.  Then at least three
of the four numbers
$$
e^{2 \pi i \alpha}, \phantom{m} E_2(\alpha), 
\phantom{m} E_4(\alpha), 
\phantom{m} E_6(\alpha)
$$
are algebraically independent.
\end{thm}
Here for $z \in \gH$, we have
$$
E_2(z) := 1 - 24 \sum_{n=1}^{\infty} \sum_{d|n} d q^n .
$$
It is a holomorphic function on $\gH$ and almost transforms like
a modular form. 
 
Furthermore, Nesterenko (\cite{NP}, page 31) suggests 
the following general conjecture which generalizes both
his and Schneider's theorem.

\begin{conj}{\rm (Nesterenko)}
Let $\alpha \in \gH$ and assume that at most three
of the following five numbers
$$
\alpha, \phantom{m} e^{2 \pi i \alpha}, 
\phantom{m} E_2(\alpha), 
\phantom{m} E_4(\alpha), \phantom{m} E_6(\alpha)
$$
are algebraically independent. Then $\alpha$ is
necessarily a  CM point and the field 
$$
{\bQ}(e^{2\pi i \alpha}, E_2(\alpha),E_4(\alpha),E_6(\alpha))
$$ 
has transcendence degree $3$.
\end{conj}

The following theorem by the first author, Murty and Rath (see \cite{GMR}) 
will play a significant role in proving Theorem \ref{three}.

\begin{thm}\label{req}
Let $\alpha \in \gH$ be such that $j(\alpha) \in \bQ$. 
Then $e^{2\pi i \alpha}$ and $\Delta(\alpha)$ are 
algebraically independent.
\end{thm}

\section{Proofs of Theorems}

\medskip

From now on, every WH modular form 
is assumed to be non-zero and with algebraic Fourier coefficients
unless otherwise stated.  

\smallskip

For any subfield ${\rm F}$ of $\bQ$ and  any  WH modular form 
$f \in M_k^{!}(\rm F)$, we define an associated function $g_f$
given by
$$
g_f (z): = \frac{ f^{12} (z)}{\Delta^k(z)}, \phantom{m} \text{ where } z \in \gH.
$$
Clearly $g_f$ is a WH modular form of weight $0$
and hence a rational function in $j$.  Since $\Delta$ does not 
vanish on $\goth H$, it follows that $g_f$ is a polynomial in $j$
with coefficients in $\rm F$. We will denote this polynomial by $P_f$.

\subsection{Proofs of Theorem \ref{one} and Theorem \ref{two}}

Let $f \in M_k^{!}(\bQ)$ be a non-zero WH modular form
such that $f \not \sim \Delta$. Hence $P_f$ is a non-constant 
polynomial with algebraic coefficients.
Now if $g_f(\alpha)$ is algebraic, we have $j(\alpha)$ is algebraic.
Moreover, if $\alpha$ is algebraic, we have by Schneider's
theorem that $\alpha$ is CM. Conversely, if $\alpha$ is
CM, then $j(\alpha)$ is algebraic and hence
$g_f(\alpha)$ is algebraic. This completes the proof
of Theorem \ref{one}.

\medskip

If $\alpha$ is an algebraic zero of $f$, from the previous argument
it follows that $j(\alpha) \in \bQ$ and hence by using 
Schneider's theorem, we get Theorem~\ref{two}.

\subsection{Proofs of Theorem \ref{three} and Theorem \ref{four}}

Let $\alpha \in \gH$ be such that $j(\alpha) \in \bQ$ and $f$ be as
in Theorem \ref{three}.
Since $j(\alpha)$ is algebraic, by Theorem \ref{req}, 
we know that $e^{2\pi i \alpha}$ and $\Delta(\alpha)$ 
are algebraically independent.
Now suppose that $f(\alpha)$ is not equal to zero. 
Since the non-zero number $g_f(\alpha)$ is a 
polynomial in $j(\alpha)$ with algebraic coefficients, it is algebraic.
Hence $e^{2\pi i \alpha}$ and $f(\alpha)$ are algebraically 
independent. This completes the proof of Theorem \ref{three}.

\medskip

We now prove Theorem \ref{four}.
Let $\alpha \in \gH$ be a non-CM algebraic number
and $f_1, f_2 \in {\rm S}_\alpha$ be WH
modular forms of weight $k_1$ and $k_2$ respectively.
By Theorem \ref{two}, neither $f_1(\alpha)$ nor 
$f_2(\alpha)$ is equal to zero. For $z \in \gH$,
consider the WH modular form  
$$
F(z):= f_1^{k_2}(\alpha) f_2^{k_1}(z) 
- f_2^{k_1}(\alpha)f_1^{k_2}(z)
$$ 
of weight $k_1k_2$. Again by Theorem \ref{two}, any zero of this WH
modular form is either CM or transcendental. Since $\alpha$ is a non-CM 
algebraic number, we get a
contradiction unless $F$ is identically zero.  Hence
$f_1 \sim f_2$. This completes the proof.

\subsection{Proofs of Theorem \ref{five} and Corollary \ref{six}}

Let $f$ be as in Theorem \ref{five}. Without loss of generality,
we can assume that $a(n_0)= 1$. Hence $P_f$ is a monic polynomial with 
algebraic coefficients.
Note that for any $\alpha \in {\goth H}$,
$$
f(\alpha)= 0  \iff  P_f(j(\alpha)) = 0.
$$
We will now prove Theorem \ref{five} by method of contradiction. 
Suppose that all zeros of $f$ are algebraic. Then 
by Theorem \ref{two} they are CM points. 
Further, we know from the theory of
complex multiplication that $j(\alpha)$ is an
algebraic integer when $\alpha$ is CM.
Since not all Fourier coefficients of $f$ are algebraic integers,
the same is true for the polynomial $P_f$. But this can not be
true if all zeros of $P_f$ are algebraic integers. This completes 
the proof of Theorem \ref{five}.

\medskip

Let $f$ be as in Corollary \ref{six}. Now if $f \in M_k^{!} (\bQ)$ 
and $f \not\in M_k^{!} (\bZ)$, then by Theorem \ref{five},
$f$ has at least one transcendental zero, a contradiction. 
This completes the proof of Corollary \ref{six}.

\subsection{Proof of Theorem \ref{seven}}

Let $f \in M_k^!(\Q)$ be a non-zero WH modular
form with all its zeros on the arc $A$.  Then
$P_f$ is a non-zero polynomial with rational coefficients.  
Let $\alpha_0 \in \gH$ be an algebraic zero of $f$. Then $j(\alpha_0)$ is
algebraic and hence $\alpha_0$ is a CM point.  Let $|D|$ be
the absolute value of the discriminant of the lattice $\Lambda := (\alpha_0 , 1)$. 
Then by Proposition \ref{cm}, there is an element 
$\sigma \in {\rm Gal}\left(\overline{\Q(\sqrt{D})} / \Q(\sqrt{D})  \right)$
such that $\sigma(j(\alpha_0)) = j(\alpha_1)$, where
$$
\alpha_1:=  \begin{cases}
      \frac{i\sqrt{|D|}}{2} &\text{if $D \equiv 0\!\!\!\pmod 4$,}\\
      \frac{ - 1+ i\sqrt{|D|}}{2} &\text{if $D \equiv 1\!\!\!\pmod 4$.}
      \end{cases}
$$
Since $P_f(x) \in \Q[x]$, applying $\sigma$ to $P_f(j(\alpha_0)) = 0$,
we have $P_f(j(\alpha_1)) = 0$. This implies that $f(\alpha_1) = 0$.
By assumption, all zeros of $f$ are on the arc $A$. Hence
$|\alpha_1| =1$ and $-1/2 \le \Re(\alpha_1) \le 0$. Hence
the only possible values of $D$ are $D= -3, -4$,
giving $\alpha_0 = \alpha_1 = i, \rho$.
This completes the proof of Theorem \ref{seven}.

\subsection{Proof of Theorem \ref{ate}}
Let
$$
g := \prod\limits_{ \gamma \in \Gamma / \Gamma_0(p)} f|\gamma.
$$
Then $g$ is a WH modular form of level one with rational 
Fourier coefficients. This is true 
because of the $q$-expansion principle
due to Deligne and Rapoport  \cite{DR}
(Theorem 3.9, p. 304) which
tells us that if an integer weight modular form $f$ has
rational Fourier coefficients at the cusp
$i \infty$, then the Fourier expansion of $f$ at all other cusps
must also have rational Fourier coefficients. Since all
zeros of $f$ lie on $B$, it follows that
all zeros of $g$ are on arc $A$. Hence by Theorem \ref{seven},
we have that all these zeros are transcendental other than $i$ and
$\rho$. As a coset decomposition of $\Gamma$ in $\Gamma_0(p)$, 
we may choose the $p+1$ elements (see page 8 of \cite{MK})
${\rm ST}^n$ ($0 \leq n <p$) and I, where 
$$
\rm S := \begin{pmatrix} 
      0 & -1 \\
      1 &  \ 0
     \end{pmatrix},
\phantom{m}
\rm T := \begin{pmatrix}
     1 & 1 \\
     0 & 1
     \end{pmatrix}
 \phantom{m}  
 \text{and }
 \phantom{m}
 \rm I := \begin{pmatrix}
     1 & 0 \\
     0 & 1
     \end{pmatrix}.
$$
Thus all zeros of $f$ are transcendental except the possible CM
zeros in $C$.

\subsection{Proofs of Theorem \ref{eisen} and  Theorem \ref{eisen1}}

Let $p = 2, 3$. By the given hypothesis, 
$E_{k,p}(z)$ is a modular form for $\Gamma_0(p)$.
Consider the level one and weight $k_p$ modular form
$$
f_p(z) := \prod_{ \gamma \in \Gamma / \Gamma_0(p)}  E_{k,p}(z) | \gamma.
$$
By the Deligne-Rapoport $q$-expansion principle, 
the modular form $f_p$ has rational Fourier coefficients.
Then $P_{f_p}$ is a non-constant polynomial with rational coefficients.
Let $\alpha_p$ be a CM zero of $g_{f_p}$. Then by Proposition
\ref{cm}, there exists an element 
$\sigma \in {\rm Gal}\left( \overline{\Q(\sqrt{D})}/ \Q(\sqrt{D})  \right)$
such that $\sigma(j(\alpha_p)) = j(\alpha_p')$, where
$$
\alpha_p' :=  \begin{cases}
      \frac{i\sqrt{|D|}}{2} &\text{if $D \equiv 0\!\!\!\pmod 4$,}\\
      \frac{-1+ i\sqrt{|D|}}{2} &\text{if $D \equiv 1\!\!\!\pmod 4$.}
      \end{cases}
$$
Since $P_{f_p}(x) \in \Q[x]$, applying $\sigma$ to $P_{f_p}(j(\alpha_p)) =0$,
we have $P_{f_p}(j(\alpha_p')) = 0$. This implies that $f_p(\alpha_p') = 0$.
But then either $E_{k,p}(\alpha_p') = 0$ or 
$E_{k,p}(\gamma\alpha_p') = 0$,
where 
$\gamma \in \{ \rm I, S, ST \}$ if $p=2$ and $\gamma \in \{ \rm I, S, ST, ST^2 \}$ 
if $p=3$. Here $\rm S, T$ are as in proof of Theorem \ref{ate}.
By the work of Miezaki, Nozaki and Shigezumi (see \cite{MNS}),
we know that all zeros of $E_{k,p}(z)$ for $p=2,3$ lie on the arc
$$
A_p :=  \left\{ z \in \gH ~|~ -1/2 \leq \Re(z) \leq 0, 
~ |z|= 1/\sqrt{p} \right\}.
$$
Hence the only possible CM zeros of 
$E_{k,2}(z)$ in $\gF_2$ can be calculated to be 
$$
\frac{i}{\sqrt{2}},  
\phantom{m} \frac{-1 + i \sqrt{7}}{4},
\phantom{m} \frac{-1 + i}{2}
$$ 
and the only possible CM zeros of $E_{k,3}(z)$ in $\gF_3$
can be calculated to be
$$
\frac{i}{\sqrt{3}},
\phantom{m} \frac{-1 + i \sqrt{11}}{6}, 
\phantom{m} \frac{-1 + i\sqrt{2}}{3}, 
\phantom{m} \frac{-3 + i \sqrt{3}}{6}.
$$

\subsection{Proofs of Theorem \ref{eight} and Theorem \ref{ten}}

Let $J = L$ or $R$ and $f_J\in M_k^!(\Q)$ be a non-zero 
WH modular form with all its zeros on the half-line $J$.
Let $\alpha_{0, J} \in \gH$ be a CM zero of $f_J$ and $|D_J|$ be the absolute value 
of the discriminant of the lattice $\Lambda := (\alpha_{0, J} , 1)$.
Then arguing as before and using 
Proposition \ref{cm}, we see that there is
a $\sigma_J \in {\rm Gal}\left(\overline{\Q(\sqrt{D_J})} / \Q(\sqrt{D_J})  \right)$
such that $\sigma(j(\alpha_{0, J})) = j(\alpha_{1, J})$, where
$$
\alpha_{1, J}:=  \begin{cases}
      \frac{i\sqrt{|D_J|}}{2} &\text{if $D_J \equiv 0\!\!\!\pmod 4$,}\\
      \frac{ - 1+ i\sqrt{|D_J|}}{2} &\text{if $D_J \equiv 1\!\!\!\pmod 4$.}
      \end{cases}
$$
It is easy to see that $\alpha_{1, J}$ is a zero of $f_J$. 
By assumption, all zeros of $f_J$ lie on $J$ . 
Hence 
$$
\Re(\alpha_{1, J}) =  \begin{cases}
          -\frac{1}{2}  & \text{ if $J=L$},\\
          0  & \text{ if $J = R$}
          \end{cases}
$$
and we have 
$$
\alpha_{1, J} = \begin{cases}
                       \frac{ - 1+ i\sqrt{|D_L|}}{2} & \text{ if $J=L$},\\
                       \frac{ i\sqrt{|D_R|}}{2} & \text{ if $J= R$},
                      \end{cases}
$$ 
i.e. $D_L \equiv 1 \!\!\pmod{4}$ and $D_R \equiv 0 \!\!\pmod{4}$.
Since $\alpha_{0, J}$ is a CM point, it is a
root of a polynomial $a_Jx^2+ b_Jx + c_J \in \Z[x]$ with $(a_J, b_J, c_J) =1, a_J > 0$
and $b_J^2 - 4a_Jc_J < 0$. Further, $D_J = b_J^2 - 4a_Jc_J$ and 
$\alpha_{0, J} = \frac{-b_J + i \sqrt{|D_J|}}{2a_J}$.
By assumption, we have 
$$
\Re(\alpha_{0, J}) = \begin{cases}
                    - \frac{1}{2}  & \text{ if $J=L$},\\
                             0  & \text{ if $J = R$}
          \end{cases}
$$
and
$$
\Im(\alpha_{0, L}) \ge \begin{cases}
                    \sqrt{3}/{2}  & \text{ if $J=L$},\\
                             1   & \text{ if $J = R$}.
          \end{cases}
$$
Hence 
$$
b_L = a_L  \text{  and } 
1 \le a_L \le \left[\sqrt\frac{|D_L|}{3} \right]
$$
and
$$
b_R = 0  \text{  and } 
1 \le a_R \le \left[\frac{\sqrt{|D_R|}}{2} \right].
$$
This completes the proofs of the theorems.

\subsection{Proofs of Remark \ref{nine} and Remark \ref{ele}}

Let $J = L$ or $R$ and $f_J  \in M_J$. Suppose that 
$f_J$ has a CM zero on the line segment $J_2$, where
$J_2 = L_2$ or $R_2$.

When $J_2 = L_2$, by Theorem \ref{eight}, zeros of $f_J$
are in the imaginary quadratic fields $K$ with discriminant
$D = -3, -7, -11, -15$.  Next, we note that it is possible
to have $a = 2$ in Theorem \ref{eight} only when $D = -15$.
This completes the proof of the Remark \ref{nine}.

Now if $J_2 = R_2$, then by Theorem \ref{ten}, we have 
$D = -4, -8, -12$. For all these $D$'s we have $a=1$. Hence the 
possible CM zeros of $f_J$ on $R_2$ are $i, ~i\sqrt{2}$ and  $i\sqrt{3}$.
This completes the proof of the Remark \ref{ele}.

\medskip
\noindent
{\bf Acknowledgments:}
The authors thank Ram Murty and Purusottam Rath for 
going through an earlier version. The authors 
would also like to thank the referee for his/her useful
suggestions.


\medskip

\begin{thebibliography}{100}

\bibitem{AKN}
T. Asai,  M. Kaneko and H. Ninomiya,  
{\em Zeros of certain modular 
functions and an application},  Comment. Math. Univ. 
St. Paul.  {\bf 46} (1997), no. 1, 93--101. 


\bibitem{CS}
C. Chowla and A. Selberg, {\em On Epstein's 
zeta-function}, J. Reine Angew. Math. 
{\bf 227} (1967), 86--110.


\bibitem{GV}
G. V. Chudnovsky, {\em  Algebraic independence of 
constants connected with the exponential and 
the elliptic functions}, Dokl. Akad. Nauk 
Ukrain. SSR Ser. A {\bf 8} (1976), 698--701. 


\bibitem{GV1}
G. V. Chudnovsky, {\em Contributions to the theory 
of transcendental numbers}, 
Mathematical Surveys and Monographs {\bf 19}, 
Amer. Math. Soc., 1984.  


\bibitem{DC}
D. Cox,  Primes of the form $x^2 + ny^2$,
{\em John Wiley and Sons}, New York, 1989.


\bibitem{DR}
P. Deligne and M. Rapoport, {\em Les sch\'emas de 
modules de courbes elliptiques},
Modular functions of one variable, II, Lecture Notes
in Math. {\bf 349}, Springer, Berlin, 1973, 143--316,.


\bibitem{DJ}
W.Duke and P. Jenkins, {\em On the zeros and coefficients 
of certain weakly holomorphic modular forms}, 
Pure Appl. Math. Q. {\bf 4} (2008), no. 4, 1327--1340. 


\bibitem{GS}
A. Ghosh and P. Sarnak, {\em Real zeros of holomorphic 
Hecke cusp forms}, J. Eur. Math. Soc. {\bf 14}
(2012), 465--487.


\bibitem{SG} 
S. Gun, {\em Transcendental zeros of certain modular forms},
Int. J. Number Theory {\bf 2} (2006), no. 4, 549--553. 


\bibitem{GMR}  
S. Gun, M. Ram Murty and P. Rath,
{\em Algebraic independence of values of modular forms}, 
Int. J. Number Theory {\bf 7} (2011), no. 4, 1065--1074.


\bibitem{KZ}
M. Kaneko and D. Zagier, {\em Supersingular $j$-invariants, 
hypergeometric series and Atkin's orthogonal polynomials}, 
AMS/IP Studies Adv. Math. {\bf 7}, 97--126 (1998).


\bibitem{NK} 
N. Kanou, {\em Transcendency of zeros of Eisenstein series},  
Proc. Japan Acad. Ser. A Math. Sci. {\bf 76} (2000), no 5, 51--54.


\bibitem{MK}
M. I. Knopp, Modular functions in Analytic Number Theory,
{\em Chelsea Publishing company}, New York, 1993.


\bibitem{WK} 
W. Kohnen, {\em Transcendence of  zeros of 
Eisenstein series and other modular functions}, 
Comment. Math. Univ. St. Pauli {\bf 52} (2003), 
no 1, 55--57. 


\bibitem{SL}
S. Lang, Elliptic Functions, {\em Addison-Wesley}, 
Reading, MA, 1973.


\bibitem{MNS}
T. Miezaki, H. Nozaki and J. Shigezumi,  {\em On the zeros of 
Eisenstein series for $\Gamma_0^*(2)$ and $\Gamma_0^*(3)$},  
J. Math. Soc. Japan  {\bf 59} (2007), no. 3, 693--706. 


\bibitem{YN} 
Y. V. Nesterenko, {\em Modular functions and transcendence},
Math. Sb. {\bf 187} (1996), no.9, 65--96.


\bibitem{NP} 
Y. V. Nesterenko and P. Philippon (Eds), 
{\em Introduction to algebraic independence theory}, 
Lecture Notes in Math. {\bf 1752}, Springer, Berlin, 2001. 


\bibitem{KR}
K. Ramachandra, {\em Some applications of Kronecker's 
limit formulas}, Ann. of Math. (2) {\bf 80} 
(1964) 104--148. 


\bibitem{RS}
F. K. C. Rankin and H. P. F. Swinnerton-Dyer, 
{\em On the zeros of Eisenstein Series}, 
Bull. London Math. Soc. {\bf 2} (1970) 169--170.

\bibitem{ZR}
Z. Rudnick, {\em On the asymptotic distribution of 
zeros of modular forms}, Int. Math. Res. Not. 
{\bf 34} (2005), 2059--2074. 

\bibitem{TS} 
T. Schneider, {\em Arithmetische Untersuchungen 
elliptischer Integrale}, 
Math Ann. {\bf 113} (1937), no. 1, 1--13.


\bibitem{JS} 
J. Shigezumi, {\em On the zeros of the 
Eisenstein series for $\Gamma_0^*(5)$ 
and $\Gamma_0^*(7)$}, Kyushu J. Math. 
{\bf 61} (2007), no. 2, 527--549.


\bibitem{HS}
R. Holowinsky and K. Soundararajan, {\em Mass 
equidistribution for Hecke eigenforms}, 
Ann. of Math. (2) {\bf 172} (2010), 
no. 2, 1517--1528. 


\end{thebibliography}
\end{document}